\documentclass[bimj,fleqn]{w-art}
\usepackage{times}
\usepackage{w-thm}
\usepackage[colorlinks,urlcolor=black,citecolor=black,linkcolor=black]{hyperref} 
\usepackage{orcidlink}
\usepackage[authoryear]{natbib}
\usepackage{amsmath,amssymb} % amsfonts
\usepackage{amsthm,amscd} % Ths usw
\usepackage{graphicx,color}
\usepackage{url}
\usepackage{booktabs} % [c]midrule
\usepackage{setspace}

\theoremstyle{plain}

\theoremstyle{definition}

\usepackage[]{graphicx}
\chardef\bslash=`\\ % p. 424, TeXbook

\DeclareMathOperator{\dcov}{dcov}

\DeclareMathOperator{\dCov}{dCov}

\DeclareMathOperator{\dCor}{dCor}

\DeclareMathOperator{\E}{E}

\DeclareMathOperator{\Prob}{P}

\DeclareMathOperator{\Normal}{\mathcal{N}}

\newcommand{\implica}{\Rightarrow}

\newcommand{\eqv}{\Leftrightarrow}
\newcommand{\flecha}{\longrightarrow}

\newcommand{\inner}[2]{\left\langle{#1},{#2}\right\rangle}
\newcommand{\norma}[1]{\left\lVert#1\right\rVert}
\newcommand{\llaves}[1]{\left\lbrace{#1}\right\rbrace}
\newcommand{\detlim}{\underset{n\to\infty}{\flecha}}

\newcommand{\distrilim}{\overset{\mathcal{D}}{\detlim}}
\newcommand{\comb}[2]{\binom{#1}{#2}}

\newcommand{\dmu}{\,\mathrm{d}\mu}

\newcommand{\amu}{a_\mu}

\newcommand{\dtheta}{\,\mathrm{d}\theta}
\newcommand{\spx}{\mathcal{X}}
\newcommand{\spy}{\mathcal{Y}}
\newcommand{\spz}{\mathcal{Z}}

\newcommand{\dx}{d_{\mathcal{X}}}
\newcommand{\dy}{d_{\mathcal{Y}}}

\newcommand{\xd}{\left(\mathcal{X},d_{\mathcal{X}}\right)}
\newcommand{\yd}{\left(\mathcal{Y},d_{\mathcal{Y}}\right)}
\newcommand{\xsd}{\left(\mathcal{X},\sqrt{d_{\mathcal{X}}}\right)}

\newcommand{\Z}{\mathbb{Z}}

\newcommand{\Zplus}{\mathbb{Z}^{+}}
\newcommand{\Nstar}{\Zplus}%mathbb{N}^{*}}

\begin{document}
\DOIsuffix{bimj.70150}
\Volume{68}
\Issue{5}
\Year{2026}
\pagespan{1}{16}
\keywords{Association measures; Distance correlation; Epistasis; Genomics of complex diseases; Schizophrenia.\\[1pc]
%\noindent\hspace*{-4.2pc} Supporting Information for this article is available from the author or on the WWW under\break \hspace*{-4pc} \underline{http://dx.doi.org/10.1022/bimj.XXXXXXX} (please delete if not applicable)
}  %%% semicolon and fullpoint added here for keyword style

\title[Testing for epistasis with distance correlation]{Testing for genetic interactions in complex disease\\with distance correlation}

%%%% AUTHORS %%%%
\author[Castro-Prado {\it{et al.}}]{Fernando Castro-Prado\orcidlink{0000-0003-3722-2013}\footnote{Corresponding author e-mail: {\sf f.castro.prado@usc.es} }\inst{1,2,3} }

\author[]{Javier Costas\orcidlink{0000-0003-0306-3990}\inst{2}}

\author[]{Dominic Edelmann\orcidlink{0000-0001-7467-6343}\inst{3}}

\author[]{\\Wenceslao Gonz\'alez-Manteiga\orcidlink{0000-0002-3555-4623}\inst{1}}

\author[]{David R. Penas\orcidlink{0000-0002-7256-3094}\inst{1,4}}

%%%% POSTAL ADDRESSES %%%%
\address[\inst{1}]{Department of Statistics, Faculty of Mathematics, University of Santiago de Compostela, Galicia, Spain.}

\address[\inst{2}]{Psychiatric Genetics Laboratory, Santiago Health Research Institute (IDIS), University Hospital, Santiago de Compostela, Galicia, Spain.}

\address[\inst{3}]{Biostatistics Department, German Cancer Research Center (DKFZ), Heidelberg, Baden-Württemberg, Germany.}

\address[\inst{4}]{Computational Biology Laboratory, Spanish National Research Council (MBG-CSIC), Pontevedra, Galicia, Spain.}

\Receiveddate{27 April 2023} \Reviseddate{25 October 2024} \Accepteddate{12 December 2024} 

\begin{abstract}
% One paragraph. Max.: 250 words. This is ~150 words.
Understanding epistasis (genetic interaction) may shed some light on the genomic basis of common diseases, including disorders of maximum interest due to their high socioeconomic burden, like schizophrenia. Distance correlation is an association measure that characterizes general statistical independence between random variables, not only the linear one. Here, we propose distance correlation as a novel tool for the detection of epistasis from case-control data of single-nucleotide polymorphisms. On the methodological side, we highlight the derivation of the explicit asymptotic null distribution of the test statistic. We show that this is the only way to obtain enough computational speed for the method to be used in practice, in a scenario where the resampling techniques found in the literature are impractical. Our simulations show satisfactory calibration of significance, as well as comparable or better power than existing methodology. We conclude with the application of our technique to a schizophrenia genetics dataset, obtaining biologically sound insights.
\end{abstract}

%\vfill

%% maketitle must follow the abstract.
\maketitle                   % Produces the title.

%\vfill

%% If there is not enough space inside the running head
%% for all authors including the title you may provide
%% the leftmark in one of the following three forms:

%% \renewcommand{\leftmark}
%% {First Author: A Short Title}

%% \renewcommand{\leftmark}
%% {First Author and Second Author: A Short Title}

%% \renewcommand{\leftmark}
%% {First Author et al.: A Short Title}

%% \tableofcontents  % Produces the table of contents.

\section{Introduction}\label{Intro}
The role of heredity in psychiatry has been studied for almost a century, with \citet{Pearson} not having ``the least hesitation'' in asserting its relevance. Today it is known that a majority of psychiatric disorders are multifactorial, complex traits. They occur as a result of a combination of genetic and environmental factors, none of which are necessary or sufficient. Furthermore, the individual effect of each of them is generally trifling. More precisely, the genome can explain up to 80 \% of the susceptibility to some of these diseases, like schizophrenia \citep{Sullivan:PGC}.

The genetic susceptibility to a psychiatric disorder is spread across a large number of variants along the genome. Although the specialized literature usually focuses simply on additive models \citep{Purcell}, biological knowledge suggests that gene-gene interactions (or \emph{epistasis}) could be one of the factors that explain the phenomenon of  \emph{missing heritability}, which contributes to the inefficiency of genome-wide association studies (GWAS) when it comes to explaining causality of complex diseases \citep{Manolio,Brandes}. Evidence from studies on model organisms also supports the importance of genetic interactions in the understanding of complex traits \citep{Mackay:Moore}.

We are interested in datasets of case-control GWASs (i.e., a collection of genotypes of ``healthy'' individuals and ``patients'') for schizophrenia. The statistical challenge hinges on using these data to detect pairs of genetic variants that significantly increase or decrease the susceptibility to develop schizophrenia, which further research can confirm with biological criteria.

These data correspond to \emph{single-nucleotide polymorphisms} (SNPs), which are variants on one of the ``letters'' of the DNA (i.e., each of them occurs at one specific point of the genome). Given that we will only consider autosomal variants, each individual can carry 0, 1 or 2 copies of the \emph{minor} allele (the least frequent of the two variants) in their diploid genome. This setting requires performing statistical inference in a context of high dimension and low sample size, where the covariates are \emph{ternary} (discrete with support of cardinality $3$).

Supporting Information Appendix A reviews the extremely wide variety of approaches to the detection of epistasis that can be found in the biostatistical literature. The main conclusion of that study effort is that there is no clear winner among the different available techniques, which justifies the maintained interest in this problem over the past few years.

On Table S1, it is shown that one conspicuous epistasis detector \citep{EpiBlaster} is based on scanning for differential behaviors of (Pearson's) correlations between cases and controls. This is unsurprising, since several authors \citep{De la Fuente,Camacho,Haeseleer} support the idea of correlation tests in this context when the data are continuous (gene expression, metabolomics, and so forth); which, however, is not the case of SNPs (ternary variables).

Moreover, such techniques usually rely on the normality of the covariates, a hypothesis that might be excessively restrictive in some cases. Therefore, the procedure by \citet{Cai:Liu} contains an interesting approach, as they manage to establish a rigorous theoretical framework for the kind of correlation tests that are convenient for epistasis detection. This technique is part of the hot topic of hypothesis testing on high-dimensional covariance structures, which has been developed almost from scratch during the past few years \citep{AnnuRev:Cai}. A short technical introduction to large-scale correlation tests (LCTs) can be found in Supporting Information Appendix C.

\begin{figure}%[!p]
	\centering\includegraphics[width=0.8\textwidth]{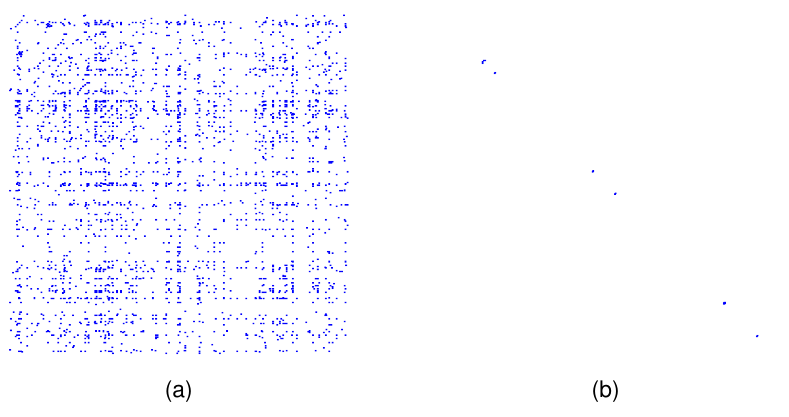}
	\caption{Adjacency matrix of the putative epistatic network detected by the LCT method encoded as ``LCT-B'', for (a) gene expression data for prostate cancer (Broad Institute) and (b) SNP data for schizophrenia (Health Research Institute, Santiago de Compostela)}
	\label{Mcoexpr_LCT-B}
\end{figure}

We have implemented the LCTs of \citet{Cai:Liu} in the \emph{R} programming language to further illustrate the motivation of our work. We first reproduced the real-data example from the original LCT article, obtaining the adjacency matrix in Figure~\ref{Mcoexpr_LCT-B}a. To accomplish that, we applied the method known as LCT-B to the data by \citet{Singh}, in which dimensionality was trimmed down to 500 using the Welch--Satterthwaite test (Behrens--Fisher problem). Since the variables involved are assumed to be continuous, the LCT-B yields believable results, in the sense that the resulting matrix is sparse, but not too much. However, a biological validation of all those results would be extremely difficult to carry out.

On the other hand, when the schizophrenia SNP data (remarkably discrete) are analyzed, the adjacency matrix looks very different (Figure~\ref{Mcoexpr_LCT-B}b) from the previous one (Figure~\ref{Mcoexpr_LCT-B}a). The only nonzero elements are close to the diagonal, owing to the fact that the only pairs that are being detected are in linkage disequilibrium (i.e., the frequency of such SNP pairs is significantly different from the product of the marginal frequencies, due to physical proximity within a certain chromosome). Such findings are useless from the point of view of psychiatric genetics because they do not show an association that is related to schizophrenia, but rather one that is independent of this disorder.

Some authors, like \citet{EpiBlaster}, argue that treating clearly discrete SNP data as continuous is an acceptable simplification. Nevertheless, even if that could be anecdotally true in some specific setting, this is clearly not the case, as Figure~\ref{Mcoexpr_LCT-B}b clearly displays.

The unsatisfactory behavior of the LCT-B, one of the most robust techniques for epistasis detection (when applied to a different setting from the one it was originally intended for), is the main motivation of the present article. In this context, it is justified to wonder which association measures characterize the independence of ternary variables, as well as how to extend the LCTs by \citet{Cai:Liu} to less stringent conditions so that they become applicable to SNP data.

The remainder of the paper is organized as follows. Section 2 summarizes the state of the art in the characterization of general statistical independence in metric spaces, focusing on an association measure named \emph{distance correlation}. Section 3 introduces a novel testing procedure for independence in ternary data. Some results of our simulation study are reported in Section 4. In Section 5, we apply the method to a genomic dataset of schizophrenia. Concluding remarks are given in Section 6.

\section{Distance correlation in metric spaces}
The \emph{energy of data} \citep{TEOD} is a branch of mathematical statistics that has been recently developed and it includes the characterization of statistical independence in Euclidean spaces via an association measure known as \emph{distance correlation}.

The extension of distance correlation to metric spaces is a nontrivial issue, which will be concisely summarised hereinafter. For a more detailed review of this theoretical framework, please refer to \citet{FW}.

\subsection{Distance correlation in Euclidean spaces}
When two random elements (vectors) $\mathbf{X}$ and $\mathbf{Y}$ are Euclidean-space-valued (let $\mathbf{X}$ be $L-$dimensional and $\mathbf{Y}$ be $M-$dimensional, for $L,M\in\Zplus$), it is possible to define an association measure that characterizes their independence with generality, which is called \emph{distance correlation} \citep{SRB}. First, distance covariance should be defined, as a weighted $L^2$ norm of the difference of the joint characteristic function and the product of the marginals:
\[\dCov(\mathbf{X},\mathbf{Y}):=\norma{\varphi_{\mathbf{X},\mathbf{Y}}-\varphi_\mathbf{X}\varphi_\mathbf{Y}}_w\equiv\sqrt{\int_{\mathbb{R}^L\times\mathbb{R}^M}|\varphi_{\mathbf{X},\mathbf{Y}}(t,s)-\varphi_\mathbf{X}(t)\varphi_\mathbf{Y}(s)|^2w(t,s)\,\mathrm{d}t\,\mathrm{d}s}\text{;}\]
where $w$ is a weight function which is dependent on the dimension of the Euclidean spaces in which the supports of $\mathbf{X}$ and $\mathbf{Y}$ are contained. As usually: \[\varphi_\mathbf{X}(t):=\E\left[\exp\left( i\inner{t}{\mathbf{X}}\right) \right],\:t\in\mathbb{R}^L\text{.}\]

Logically, distance correlation is defined as the quotient of covariance and the product of standard deviations (as long as none of the latter vanish):
\[\dCor(\mathbf{X},\mathbf{Y}):=\frac{\dCov(\mathbf{X},\mathbf{Y})}{\sqrt{\dCov(\mathbf{X},\mathbf{X})\dCov(\mathbf{Y},\mathbf{Y})}}\text{,}\]
and so it has no sign. It is an improved version of Pearson's correlation squared because it has values in [0,1] and, more importantly, it is zero if \emph{and only if} $\mathbf{X}$ and $\mathbf{Y}$ are (statistically) independent.

\begin{sloppypar}
	However convoluted the initial definition of \emph{dCor} is, the sample version can easily be computed. Given a paired sample $(\mathbf{X}_1,\mathbf{Y}_1),\ldots,(\mathbf{X}_n,\mathbf{Y}_n)\text{ IID }(\mathbf{X},\mathbf{Y})\text{;}$ let $a_{ij}:=d(\mathbf{X}_i,\mathbf{X}_j)$ be the Euclidean distance between the $\mathbf{X}$'s with indices \mbox{$i,j\in\{1,\ldots,n\}$}. Then, the doubly centered distances are:
	\begin{equation} A_{ij}:=a_{ij}-\frac{1}{n}\sum_{k=1}^na_{ik}-\frac{1}{n}\sum_{k=1}^na_{kj}+\frac{1}{n^2}\sum_{k,l=1}^na_{kl}\,. \label{eq:center} \end{equation}
	If $\{b_{ij}\}_{i,j}$ and $\{B_{ij}\}_{i,j}$ are analogously defined for $\{\mathbf{Y}_i\}_i$, the empirical distance covariance is the nonnegative square root of
	\begin{equation} \widehat{\dCov}_n(\mathbf{X},\mathbf{Y})^2:=\frac{1}{n^2}\sum_{i,j=1}^n A_{ij}B_{ij}\,. \label{eq:empdc} \end{equation}
\end{sloppypar}

The above estimator is reminiscent of the following alternative representation of $\dCov^2$:
\[\dCov(\mathbf{X},\mathbf{Y})^2\!=\!\E\Big[\Big(d(\mathbf{X},\mathbf{X}')-\E\lbrace d(\mathbf{X},\mathbf{X}'')\rbrace-\E\lbrace d(\mathbf{X}',\mathbf{X}''')\rbrace+\E\lbrace d(\mathbf{X}'',\mathbf{X}''')\rbrace\Big)\]
\[\times\Big(d(\mathbf{Y},\mathbf{Y}')-\E\lbrace d(\mathbf{Y},\mathbf{Y}'''')\rbrace-\E\lbrace d(\mathbf{Y}',\mathbf{Y}''''')\rbrace+\E\lbrace d(\mathbf{Y}'''',\mathbf{Y}''''')\rbrace\Big)\Big]\text{,}\]
which is valid as long as moments of order $2$ are finite \citep[][Remark 4.6]{Jakobsen}. In the above equation, primed letters refer to independent and identically distributed (IID) copies of the corresponding random vector, and $d$ denotes the Euclidean distance

Whenever $\{\mathbf{X},\mathbf{Y}\}$ are independent and have finite first moments, the asymptotic null distribution of a scaled version of the preceding statistic is a linear combination of independent chi-squared variables with one degree of freedom. More precisely,
\[n\:\widehat{\dCov}_n(\mathbf{X},\mathbf{Y})^2\distrilim\sum_{j=1}^\infty \lambda_j Z_j^2\text{,}\] 
where $\{Z_j\}_j$ are IID $\Normal(0,1)$ and where $\{\lambda_j\}_j\subset\mathbb{R}^+$.

Unfortunately, knowing the form of the theoretical null distribution is often not helpful in practice. As a result, all the distance correlation literature we are aware of resorts to resampling techniques when it comes to approximating the critical values for the independence test. They generally design the resampling scheme based on the information that $H_0$ provides, which in this case (i.e., independence) leads to permutation testing. The theoretical quadratic form is not used, due to the difficulty of estimating the $\lambda_j$'s, since one would need to deal with an abstract linear operator and obtain its nonzero eigenvalues (potentially, an infinity of them), all this under no model assumptions (note that our setting is nonparametric); see, for example, \citet{Jakobsen} or \citet{SRB}.

\subsection{The generalized distance covariance}
Let $\theta$ be a (Borel) probability distribution on $\spx\times\spy$, where $\xd$ and $\yd$ are separable metric spaces. Its distance covariance is defined as
\[\dcov(\theta):=\int_{(\spx\times\spy)^2}d_\mu(x,x')d_\nu(y,y')\dtheta^2\left((x,y),(x',y')\right)\text{,}\]
where $(\mu,\nu)$ are the marginals of $\theta$ on $(\spx,\spy)$ and are assumed to both have finite first moments \citep{Lyons,Jakobsen}. Functions $d_\mu$ and $d_\nu$ are the doubly centered versions of $\dx$ and $\dy$, respectively:
\begin{align*}
	d_\mu:\;&\spx\times\spx\,\longrightarrow\mathbb{R}\\
	&(x_1,x_2)\longmapsto\dx(x_1,x_2)-\amu(x_1)-\amu(x_2)+D(\mu)\text{;}
\end{align*}
where:
\begin{align*}
	\amu:\;&\spx\longrightarrow\mathbb{R}\\
	&x\longmapsto\int_{\spx}\dx(x,x')\dmu(x')\text{;}
\end{align*}
\[D(\mu):=\int_{\spx}\amu\dmu=\int_{\spx^2}\dx\dmu\times\mu\text{.}\]

Like ordinary covariance, \emph{dcov} vanishes under independence as a result of Fubini's theorem:
\[\dcov(\mu\times\nu)=\{D(\mu)-2D(\mu)+D(\mu)\}\{D(\nu)-2D(\nu)+D(\nu)\}=0\text{.}\]

\subsection{Distance covariance in negative type spaces}
The fact that
\[\theta=\mu\times\nu\implica\dcov(\theta)=0\]
makes it natural to wonder which spaces ensure that the reciprocal implication also holds. The answer is \emph{strong negative type} spaces, since in them $\dcov(\theta)$ can be presented as an injective function of $\theta-\mu\times\nu$.

A metric space $\xd$ is said to be of negative type if and only if:
\[\forall n\in\Z^+;\:\forall x,y\in\spx^n:\:2\sum_{i,j=1}^n\dx(x_i,y_j)\geq\sum_{i,j=1}^n\{\dx(x_i,x_j)+\dx(y_i,y_j)\}\text{.}\]

There are many familiar examples of negative type spaces, including all Hilbert spaces (and therefore, all Euclidean geometries).

If $\xd$ has negative type, the inequality $D(\mu_1-\mu_2)\leq0$ holds for any probability distributions $\mu_1,\mu_2$ on $\spx$ with finite first moments. On top of that, if the operator $D$ separates probability measures (with finite first moments) in $\xd$, that space is said to have \emph{strong} negative type:
\[D(\mu_1-\mu_2)=0\eqv\mu_1=\mu_2\text{.}\]
Whenever $\spx$ and $\spy$ have strong negative type,
\[\dcov(X,Y)\overset{\text{def.}}{=}\dcov(\theta)=0\eqv X,Y\text{ independent,}\]
for any random element $(X,Y)\sim\theta$ with values in $\spx\times\spy$.

\subsection{Nonparametric test of independence in metric spaces}
For $n\in\mathbb{Z}^{+}$, the \emph{empirical measure} associated with a certain sample
\[\llaves{(X_i,Y_i)}_{i=1}^n\text{ IID }(X,Y)\sim\theta\]
is defined as customarily:
\[\theta_n:=\frac{1}{n}\sum_{i=1}^n\delta_{(X_i,Y_i)}\:\text{,}\]
where $\delta_z$ denotes point mass at $z\in\spx\times\spy$.

It is easy to see that the natural estimator $\dcov(\theta_n)$ is a $V-$statistic
\[\dcov(\theta_n)=\frac{1}{n^6}\sum_{i_1=1}^n\cdots\sum_{i_6=1}^n h\left((X_{i_\lambda},Y_{i_\lambda})_{\lambda=1}^6\right)\text{,}\]
whose (nonsymmetric) kernel $h$ is given by:
\begin{alignat*}{3}
	h:\;&(\spx\times\spy)^6&\longrightarrow&\mathbb{R}\\
	&\big((x_i,y_i)\big)_{i=1}^6&\longmapsto&\left\lbrace\dx(x_1,x_2)+\dx(x_3,x_4)-\dx(x_1,x_3)-\dx(x_2,x_4)\right\rbrace\\
	&&&\times\left\lbrace\dy(y_1,y_2)+\dy(y_5,y_6)-\dy(y_1,y_5)-\dy(y_2,y_6)\right\rbrace\text{.}
\end{alignat*}

If $\theta$ is the product of its marginals and these are nondegenerate, the asymptotic null distribution of the $V-$statistic is
\[n\dcov(\theta_n)\distrilim\sum_{i=1}^\infty\lambda_i(Z_i^2-1)+D(\mu)D(\nu)\text{;}\]
where $\llaves{Z_i}_{i\in\Nstar}\text{ IID }\Normal(0,1)$ and where $\llaves{\lambda_i}_{i\in\Nstar}\subset\mathbb{R}$ are unknown (dependent on $\theta$). The most logical approach is, once again, resorting to permutation tests.

\section{Distance correlation-based test for epistasis}\label{ourtest}
Thus far, the theoretical basis for the usage of distance correlation within certain metric spaces has been set. Upon this, extensions and applications can be formulated. There are a few examples of this in the literature, in contexts such as time series \citep{Davis}, topological data analysis \citep{TDA}, and discretized stochastic processes \citep{Dehling}.

What is about to be presented is the particularization of the theoretical framework to spaces of cardinality $3$, to then design a procedure that is adapted to epistasis detection and that at the same time solves the inadequacy of the tests by \citet{Cai:Liu} to this setting (as shown in Figure~\ref{Mcoexpr_LCT-B}).

\subsection{Distance correlation in spaces of cardinality $3$}\label{sec:space3}
Clearly, in a finite space, the finiteness of moments (of any order) and separability are not an issue. In addition, one can resort to brute force and solve the system of inequations that are derived from simply using the definitions \citep{Klebanov,Lyons}, obtaining a direct \textemdash albeit cumbersome\textemdash~proof of the fact that any 3-point space $\xd$ is necessarily of strong negative type. Such proof is, in principle, superfluous, as long as one wants to make use of strong theorems, such as Schoenberg's: $\xsd$ can clearly be embedded into a Hilbert space, isometric to the vertices of a triangle in $\mathbb{R}^2$ (note that the square root transformation preserves the triangle inequality of the metric). Nevertheless, it is interesting to check that, when the metric structure becomes so simple, abstract arguments (such as the ones that arise in the proof of Schoenberg's theorem) become unnecessary.

Let $\spx:=\{0,1,2\}$ be the set of the three possible genotypes for each SNP. There is no biological reason to assume that $2\in\spx$ copies of the minor allele affect twice as much as one \citep{Bush}, neither when it comes to increasing the susceptibility to a psychiatric disorder nor to decreasing it. As a matter of fact, in some cases this susceptibility is maximal under heterozygosis \citep{Costas:Heteroz_opt}, which is coded by $1\in\spx$.

Therefore, there is no rationale for prioritizing the Euclidean distance:
\[d(0,2)=2d(0,1)=2d(1,2)\]
over more general (non-``linear'') metric spaces. And this is why distance correlation turns out to be a sensible way to adapt the ideas of \citet{Cai:Liu} to our setting. As previously commented, the marked discreteness of SNP data provides another incentive for transcending the idea of linear correlation.

No specific type of interaction is being looked for---the aim of this article is to simply detect epistasis. For this reason, we will henceforward focus on the \emph{equilateral distance} (also known as the discrete metric):
\[d(0,1)=d(1,2)=d(0,2)=1\text{.}\]
Three further conspicuous distances are to be defined, which are the degenerate ones (in which two of the vertices of the triangle are the same point), that can be especially illustrative in the analysis of real data because they allow a straightforward interpretation of the allelic model that is being studied:
\begin{enumerate}
	\item Recessive (distance ``0=1''): $\:\qquad d(0,1)=0;\;\;d(0,2)=d(1,2)=1$.
	\item Heterozygous (distance ``0=2''): $\:\; d(0,2)=0;\;\;d(0,1)=d(2,1)=1$.
	\item Dominant (distance ``1=2''): $\qquad d(1,2)=0;\;\;d(1,0)=d(2,0)=1$.
\end{enumerate}

Technically, the three previous distances are not metrics, but rather pseudometrics, as they yield $d(x,y)=0$ for $x\neq y$. However, one can modify the labels on the initial 3-point space so that it becomes a 2-point space, in which case the former pseudometric is an actual metric for which all the theoretical results of \citet{Jakobsen} immediately follow. More importantly, distance covariance could also be defined for semimetric spaces (i.e., by dropping the triangle inequality), using the definition by \citet{Sejdinovic}. Therefore, having a premetric space would suffice for a notion of generalized distance covariance to exist.

\subsection{Construction of a hypothesis test}\label{sec:test}
Searching for epistasis consists in looking for significantly different dependence structures between the case and control groups, as previously discussed. To simplify notation, let $Z_i$ and $Z_j$ be random variables with support $\spz\in\llaves{\spx,\spy}$, corresponding to two different SNPs, for which a joint sample of size $n\in\mathbb{Z}^{+}$ is available:
\[(Z_{i,1},Z_{j,1}),\ldots,(Z_{i,n},Z_{j,n})\text{ IID }(Z_i,Z_j)\text{.}\]
The aim is to test the independence of $\{Z_i,Z_j\}$ or, equivalently,
\[\begin{cases}{H_{0ij}}:\dcov(Z_i,Z_j)=0\\H_{1ij}:\dcov(Z_i,Z_j)\neq0\;;\end{cases}\]
which extends the large-scale multiple tests of local covariance and correlation structures (see Supporting Information Appendix C.1) to the use of distance covariance and correlation, instead of their product-moment counterparts.

In order to approximate the null distribution of the test statistic $\widehat{\dcov}(Z_i,Z_j)$, one can take advantage of the beauty of finite marginal spaces, in which only a few of the coefficients of the quadratic form that gives the asymptotic null distribution of distance covariance will be non-null. Namely, we present two theorems for such distributions, both for the metric of maximum interest to us and for the Euclidean distance (i.e., for classical distance covariance). Proofs can be found in Supporting Information Appendix B.

\begin{theorem} \label{th:discrete}
	Let $(X_1,\ldots,X_n)$ and $(Y_1,\ldots,Y_n)$ be IID samples of jointly distributed random variables $(X,Y) \in \{0,1,2\} \times \{0,1,2\}$, with $p_j:= \Prob(X=j)$ and $q_j:=\Prob(Y=j)$, for $j\in\{0,1,2\}.$
	
	Consider $\mathcal X = \mathcal Y =\{0,1,2\}$ equipped with the discrete metric.
	
	Then, whenever $X$ and $Y$ are independent, for $n \to \infty$,
	\[
	n \, \widehat{\dCov}_{\text{discrete}}^2 \distrilim  
	\lambda_1 \, \mu_1 Z_{11}^2 + \lambda_1 \mu_2 Z_{12}^2 + \lambda_2 \mu_1 Z_{21}^2 + \lambda_2 \mu_2 Z_{22}^2;
	\]
	
	where $Z_{11}^2, Z_{12}^2, Z_{21}^2,$ and $Z_{22}^2$ are independently chi-squared distributed with one degree of freedom. $\lambda_1$ and $\lambda_2$ are given by
	\begin{equation}\label{eq:lambda}
	\lambda_{1,2} = \frac{1-\sum p_j^2}{2} \pm \sqrt{\frac{(1-\sum p_j^2)^2}{4} - 3 \prod p_j}. %\quad  \lambda_2 = \frac{1-\sum p_j^2}{2} - \sqrt{\frac{(1-\sum p_j^2)^2}{4} - 3 \prod p_j}.
	\end{equation}
	Similarly $\mu_1$ and $\mu_2$ are given by
	\begin{equation}\label{eq:mu}
	\mu_{1,2} = \frac{1-\sum q_j^2}{2} \pm \sqrt{\frac{(1-\sum q_j^2)^2}{4} - 3 \prod q_j}. %\quad \quad \mu_2 = \frac{1-\sum q_j^2}{2} - \sqrt{\frac{(1-\sum q_j^2)^2}{4} - 3 \prod q_j}.
	\end{equation}\qed
\end{theorem}

\begin{theorem}\label{th:Euclidean}
	Let $(X_1,\ldots,X_n)$ and $(Y_1,\ldots,Y_n)$ be IID samples of jointly distributed random variables $(X,Y) \in \{0,1,2\} \times \{0,1,2\}$, with $p_j:= \Prob(X=j)$ and $q_j:=\Prob(Y=j)$, for $j\in\{0,1,2\}.$
	
	Consider $\mathcal X = \mathcal Y =\{0,1,2\}$ equipped with the Euclidean metric.
	
	Then, whenever $X$ and $Y$ are independent, for $n \to \infty$,
	\[
	n \, \widehat{\dCov}_{\text{Euclidean}}^2 \distrilim  
	\lambda_1 \, \mu_1 Z_{11}^2 + \lambda_1 \mu_2 Z_{12}^2 + \lambda_2 \mu_1 Z_{21}^2 + \lambda_2 \mu_2 Z_{22}^2;
	\]
	
	where $Z_{11}^2, Z_{12}^2, Z_{21}^2,$ and $Z_{22}^2$ are independently chi-squared distributed with one degree of freedom. $\lambda_1$ and $\lambda_2$ are given by
	\[
	\lambda_{1,2} = p_0(1-p_0)+p_2(1-p_2) \pm \sqrt{\Big(p_0(1-p_0)+p_2(1-p_2)\Big)^2 - 4 \prod p_j}. %\quad  \lambda_2 = \frac{1-\sum p_j^2}{2} - \sqrt{\frac{(1-\sum p_j^2)^2}{4} - 3 \prod p_j}.
	\]
	Similarly $\mu_1$ and $\mu_2$ are given by
	\[
	\mu_{1,2} = q_0(1-q_0)+q_2(1-q_2) \pm \sqrt{\Big(q_0(1-q_0)+q_2(1-q_2)\Big)^2 - 4 \prod q_j}. %\quad \quad \mu_2 = \frac{1-\sum q_j^2}{2} - \sqrt{\frac{(1-\sum q_j^2)^2}{4} - 3 \prod q_j}.
	\]\qed
\end{theorem}

We will now show that Theorems \ref{th:discrete} and \ref{th:Euclidean} induce asymptotic tests for testing the null hypothesis that $X$ and $Y$ are independent.

\begin{corollary}\label{cor:test:plug}
	Let $(X_1,\ldots,X_n)$ and $(Y_1,\ldots,Y_n)$ be IID samples of jointly distributed random variables $(X,Y) \in \{0,1,2\} \times \{0,1,2\}$. For each $j\in\{0,1,2\}$, we define \[\widehat{p}_j := \frac{1}{n} \sum_{i=1}^n 1_{\{X_i = j\}};\; \widehat{q}_j =  \frac{1}{n} \sum_{i=1}^n 1_{\{Y_i = j\}},\] and let $\widehat{\lambda}_{1,2}$ and $\widehat{\mu}_{1,2}$ be the plug-in versions of Equations \eqref{eq:lambda} and \eqref{eq:mu} by $\widehat{p}_j$ and $\widehat{q}_j$, respectively, that is,
	\[
	\widehat{\lambda}_{1,2} = \frac{1-\sum \widehat{p}_j^2}{2} \pm \sqrt{\frac{(1-\sum \widehat{p}_j^2)^2}{4} - 3 \prod \widehat{p}_j};
	\]
	\[
	\widehat{\mu}_{1,2} = \frac{1-\sum \widehat{q}_j^2}{2} \pm \sqrt{\frac{(1-\sum \widehat{q}_j^2)^2}{4} - 3 \prod \widehat{q}_j}.
	\]
	Denoting the joint sample by $\mathcal{S} = \{(X_1,Y_1),\ldots,(X_n,Y_n)\}$, given a nominal level $\alpha\in]0,1[$, we define the threshold $\widehat{c}_\alpha(\mathcal{S})$ as the largest nonnegative number satisfying:
	\[
	\Prob(\widehat{\lambda_1} \, \widehat{\mu}_1 Z_{11}^2 + \widehat{\lambda_1}  \widehat{\mu}_2 Z_{12}^2 + \widehat{\lambda_2} \widehat{\mu_1} Z_{21}^2 + \widehat{\lambda}_2 \widehat{\mu}_2 Z_{22}^2 \geq \widehat{c}_\alpha(\mathcal{S}) \ | \mathcal{S}) = \alpha,
	\]
	where we are simplifying the notation by writing $\widehat{\lambda_1} \equiv \widehat{\lambda_1} (\mathcal{S})$ and so forth.
	Then, if $X$ and $Y$ are independent, it holds that
	\[
	\Prob(n \, \widehat{\dCov}_{\text{discrete}}^2(\mathcal{S}) \geq  \widehat{c}_\alpha(\mathcal{S}))  \longrightarrow \alpha,
	\]
	as $n \to \infty$.
	Hence the function \[\varphi: \mathcal{S}\mapsto\Prob\left\lbrace  \widehat{\dCov}_{\text{discrete}}^2(\mathcal{S}) \geq  \widehat{c}_\alpha(\mathcal{S}) \right\rbrace \] is an asymptotic test of size $\alpha$ for testing the nonparametric null hypothesis
	\[
	H_0: \text{$X$ and $Y$ are independent}
	\]
	against the alternative
	\[
	H_1: \text{$X$ and $Y$ are not independent.}
	\]\qed
\end{corollary}

The only remaining task is to calculate the quantiles $\widehat{c}_\alpha$, for which efficient algorithms are available \citep{Duchesne}.

An analogous result to Corollary~\ref{cor:test:plug} holds for the test based on $n \, \widehat{\dCov}_{\text{Euclidean}}^2$, which would also work by plugging in the allele frequencies into the expression of the coefficients of the Gaussian quadratic form.

It is crucial to note that we do not want to directly test for the equality of distance correlations, as \citet{Cai:Liu} did when looking for differential gene co-expression, following the rationale by \citet{De la Fuente} and others. In our search for epistasis, however, we are just interested in finding SNP pairs that are dependent for the cases and independent for the controls. When, for some SNP pair, independence is rejected for healthy individuals and not for patients, it will be attributed to a spurious interaction resulting from \emph{population substructure} \citep{Brandes}, that is, from the effect of unmeasured (and often unmeasurable) covariates.

\section{Simulations}\label{sec:simu}

In order to validate our testing procedure, we have designed some population models in which the intensity of dependence can be adjusted by tuning a parameter. We first introduce those models, to then use them to compare the performance of our method with that of BOOST \citep{BOOST}, one of the most popular epistasis detectors within the genomics community. For distance covariance, we will consider the discrete metric in every scenario because it reflects our agnosticism on the underlying genetic model.

\subsection{Design of population models for the validation of our methodology}
The theoretical models that are about to be defined refer to the interaction between an arbitrary pair $\{Z_i,Z_j\}$, where $Z$ is either $X$ or $Y$, depending on the case. When it came to setting the marginal frequencies, instead of allowing for two degrees of freedom on each marginal, a further restriction was introduced (apart from the sum being one): allele and genotype frequencies were constrained to be in Hardy--Weinberg equilibrium \citep{Hardy}, since all the SNPs in the schizophrenia database verify it (it is one of the quality controls that are used). Hence, there is a single free parameter, which is the minor allele frequency, that is sampled from a uniform distribution on $[0.05,0.2]$. The lower limit mimics standard GWAS quality control filters (in settings with moderate sample size) and the upper one was set so that the resulting true interactions are not the easiest to detect.

\setcitestyle{square}There are a few options in the literature for simulating epistasis between SNPs. Some models (like the ones by \citet{Marchini}) are overly simplistic, for example by not allowing to adjust the interaction intensity in order to assess the robustness against different alternatives. Some recent approaches (like the ones studied by \citet{Russ})\setcitestyle{round} make interpretability more difficult, in the sense that we are very interested in quantifying the intensity of interaction (i.e., deviation from the null hypothesis) when assessing the power of our test. In order to overcome such shortcomings, we introduce our own models for SNP-SNP interaction.

The most straightforward model is one in which the probability of each genotype is the product of the marginals (there is independence). For dependence, two kinds of models will be defined. On the one hand, model \texttt{qexp} conveys a dependence structure that becomes more intense as parameter $e\in[1,+\infty[$ increases, in the way that Table~\ref{qexp} describes. On the other hand, model \texttt{qmult} has $g\in[0,1]$ as its free parameter (Table~\ref{qmult}). Again, the closer the parameter is to $1$, the less notorious the association becomes.

\begin{table}%
	\centering
		\caption{Contingency table for model \texttt{qexp}}
		\begin{tabular}{cccc|c}
			\multicolumn{1}{c|}{$Z_i$ \textbackslash $\:Z_j$} & 0     & 1     & 2     &  \\
			\cmidrule{1-4}    \multicolumn{1}{c|}{0} & $pr+q^es-qs$ & $ps-q^es+qs$ & $p(1-r-s)$ & $p$ \\
			\multicolumn{1}{c|}{1} & $qr-q^es+qs$ & $q^es$ & $q(1-r-s)$ & $q$ \\
			\multicolumn{1}{c|}{2} & $(1-p-q)r$ & $(1-p-q)s$ & $(1-p-q)(1-r-s)$ & $1-p-q$ \\
			\midrule
			& $r$   & $s$   & $1-r-s$ & $1$ \\
		\end{tabular}%
	\label{qexp}
\end{table}

\begin{table}%
	\centering
		\caption{Contingency table for model \texttt{qmult}}
		\begin{tabular}{cccc|c}
			\multicolumn{1}{c|}{$Z_i$ \textbackslash $\:Z_j$} & 0     & 1     & 2     &  \\
			\cmidrule{1-4}    \multicolumn{1}{c|}{0} & $pr-(1-g)qs$ & $ps+(1-g)qs$ & $p(1-r-s)$ & $p$ \\
			\multicolumn{1}{c|}{1} & $qr+(1-g)qs$ & $gqs$ & $q(1-r-s)$ & $q$ \\
			\multicolumn{1}{c|}{2} & $(1-p-q)r$ & $(1-p-q)s$ & $(1-p-q)(1-r-s)$ & $1-p-q$ \\
			\midrule
			& $r$   & $s$   & $1-r-s$ & $1$ \\
		\end{tabular}%
	\label{qmult}%
\end{table}%

\subsection{Results of the simulation study}
Each simulation consisted in the study of one of the models for a SNP pair. This is an acceptable simplification because the current setting is a problem of multiple testing and not a single high-dimensional test (see Supporting Information Appendix C.1). This means that it is not necessary to build and replicate a whole $n\times L$ matrix, and that it suffices to perform replications at the SNP-pair level.

\begin{figure}%[!p]
	\centering\includegraphics[width=.6\textwidth]{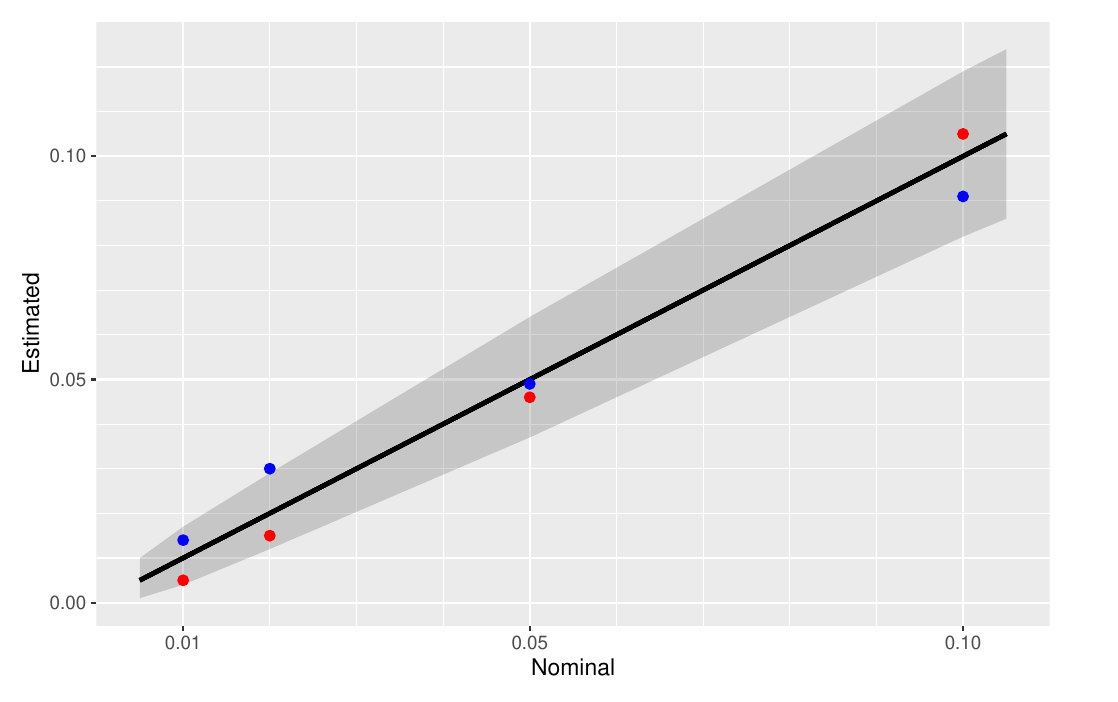}
	\caption{Empirical power under the null hypothesis ($\hat\alpha$) versus nominal significance level ($\alpha$), under model \texttt{indep} in cases and \texttt{indep} in controls. Blue dots correspond to \textit{dcov}; the red ones were generated with BOOST. The gray shadow is a 95 \% confidence band for $\hat\alpha$ given $\alpha$.}
	\label{Calibr}
\end{figure}

\begin{figure}%[!p]
	\centering\includegraphics[width=\textwidth]{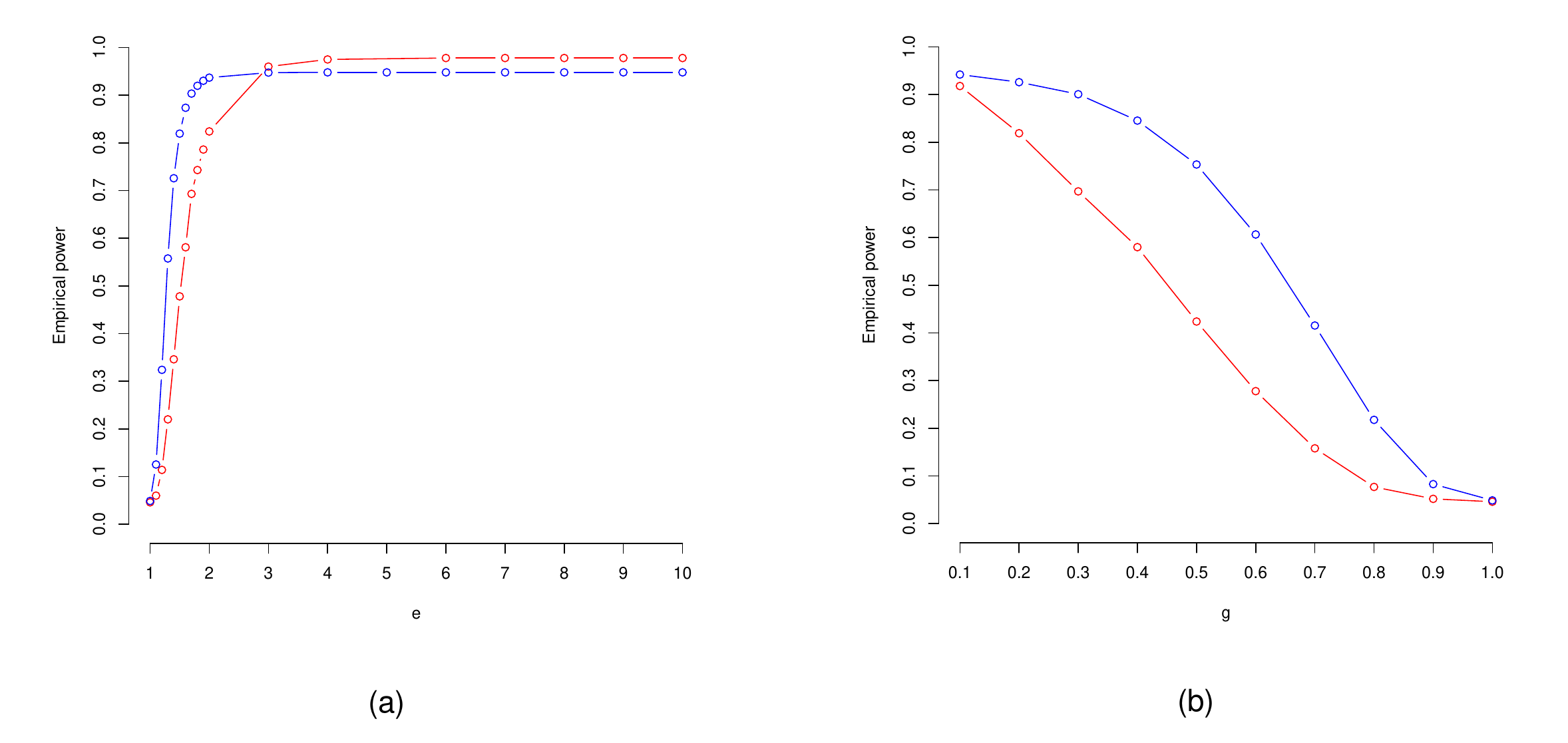}
	\caption{Empirical power when the SNP pair distributions for cases/controls are (a) \texttt{qexp} with parameter $e\in\mathbb{Z}^{+}$ and \texttt{indep}, and (b) \texttt{qmult} with parameter $g\in[0,1]$ and \texttt{indep}. Color blue represents our distance covariance test; whereas red corresponds to BOOST.}
	\label{Power_indep_q}
\end{figure}

We now briefly show some simple, illustrative examples of the performance of our testing procedure. First, Figure~\ref{Calibr} shows the calibration of significance for some usual nominal levels for the only scenario under the null hypothesis we expect to come across in practice, that is, independence in both cases and controls. On the other hand, empirical power is represented in Figure~\ref{Power_indep_q}. In all cases, $R=1000$ replicates were carried out. The sample size in each replicate, for each of the control and case group, is $500$ plus a perturbation that is uniformly drawn from the integer interval $[-10,10]\cap\Z$, in order to demonstrate the applicability of our method in scenarios where $n_1\neq n_2$.

In each plot, we also display the results we obtained with one of the most popular tools within the genomics community for the kind of epistasis we are studying --- it is called \textit{BOOST} \citep{BOOST} and is easily accessible from the widely used genetics software package PLINK \citep{PLINK}. As indicated in Supporting Information Appendix A, there is an extremely large number of options in the literature to perform this task and therefore it is not feasible to compare our technique with a representative fraction of them.

On the basis of Figure~\ref{Calibr}, it can be concluded that the calibration of significance is acceptable or even good for the most usual levels of nominal $\alpha$. In addition, the plots in Figure~\ref{Power_indep_q} show that the power is very satisfactory and that, as expected, it increases as one gets further away from the null hypothesis. In the scenarios we studied, we have either comparable or more power than BOOST.

\section{Application to a case-control study of schizophrenia}\label{SCZ}

The genomic database that we study in this article is described in detail in Supporting Information Appendix F, for the sake of reproducibility. It contains observations of $6,371,078$ SNPs across the entire genome, from a case-control study of schizophrenia in Galicia \citep{Galicia}, with $n_1=585$ cases and $n_2=573$ controls. For a better understanding of the nature of the dataset and the quality controls \citep{gwas:bimj} and downstream analysis it underwent, we refer the reader to Supporting Information Appendices F and G, the latter of which also contains further details on reproducibility.

In Section~\ref{experim0}, we illustrate the testing procedure for genetic interaction for which we set the theoretical basis in Section~\ref{ourtest} and whose numerical performance we examined in Section~\ref{sec:simu}. Then, in Sections~\ref{experim1} and \ref{experim2}, we present two experimental setups we carried out to better understand the performance of our methodology in biostatistical practice. In them, we use our DNA data as the starting point for interpretations at ``higher'' levels of the biomolecular hierarchy (proteins and RNA), based on missense SNPs and genetically regulated gene expression, respectively. As with the simulations, we restrict ourselves to the discrete metric, in order to be agnostic regarding the underlying genetic model. In the latter two experiments, we apply our methodology to SNP pairs across the human genome, to then interpret the results by performing several tests comparing proportions and ranks in cases versus controls. The key rationale is that the set of putative interactions detected with our testing procedure (i.e., the SNP pairs for which independence is rejected in cases and not in controls) will include both pairs in ``true'' epistasis and noise, whereas the SNP pairs where independence is rejected in controls and not in cases are just spurious interactions (due to population substructure).

\subsection{Illustration of the distance correlation-based test for genetic interaction}\label{experim0}

As an example of how our testing procedure would work out when testing for epistasis between two SNPs, let us consider SNP {rs6972158} (in chromosome 7) and SNP {rs61749697} (chromosome 20). These are missense SNPs (i.e., the variation of that nucleotide changes the aminoacid sequence of the protein encoded by the gene) in genes \textit{NPSR1} and \textit{GNAS}, respectively. \textit{NPSR1} codes for a G-protein-coupled receptor for neuropeptide S. \textit{GNAS} codes for a protein that binds to many different G-protein-coupled receptors, being essential for the transduction of cellular signals. Using our schizophrenia dataset, and applying the distance covariance test, we get a $p$-value of 0.041 in cases and of 0.966 in controls. This means that, for nominal $\alpha$ of 0.05, we reject independence in cases and not in controls. Then this pair would be in putative epistasis.

In order to check if this makes sense biologically, we conducted a search in the STRING database version 12.0 \citep{STRING}. According to the interaction prediction pipeline of STRING, based on several sources of evidence, the combined confidence score of the functional interaction between these two proteins is 0.619, regarded as a medium level. Confidence scores are scaled between zero and one, and correspond to the estimated likelihood of the association under consideration being true, given the underlying evidence. Therefore, the putative interaction between these two SNPs detected by our method is biologically plausible.

\subsection{Experiment I: Functional enrichment}\label{experim1}

Taking into account the goal of this first experiment, for the reasons stated in Supporting Information Appendix G, it is sensible to restrict ourselves to a certain subset of the initial database, comprising $L=8030$ missense SNPs.

First, we apply our test procedure separately to cases and controls (as previously discussed), with a \citet{BH} nominal threshold for the false discovery rate (FDR) of 0.05. We only consider SNP pairs consisting of two variants that lie on different chromosomes or that are more than 1 Mb apart (i.e., not physically close). This prevents evident cases of spurious findings due to linkage disequilibrium \citep{BOOST}.

We thus obtain 113 out of $\comb{L}{2}$ SNP pairs that show association in cases and not in controls (which we consider putative interactions), versus 95 in controls and not in cases (which just reflect population substructure). The difference (in proportions) is not significant; with a $p$-value of $0.12$, which is greater than any conventionally used nominal value of $\alpha$. These 113 and 95 pairs correspond, respectively, to 222 and 189 unique SNPs, a proportion difference with $p$ of 0.055. Those SNPs lie on 220 and 191 different genes. Removing the 13 that are common among both lists, we get 207 and 178 genes ($p$ of 0.07).

We hypothesize that there will be an overrepresentation of genes involved in synapse, which is a brain function known to be closely related to schizophrenia \citep{SynGO}, in our group of putative interactions with respect to the spurious ones. Intersecting, we see that 13 of the 207 and 9 out of 178 genes are known to be related to synapse. The proportion difference has a $p$-value of 0.26, so our results for this part are negative --- we can show no (strong) evidence that we are detecting any signal related to synapse.

However, this is not to say that our method cannot offer interesting insight into these data. The current knowledge on complex disease genetics indicates that regulatory regions play a crucial role \citep{Sullivan:Geschwind}, so one should focus on genetically regulated gene expression, rather than on missense SNPs. This motivates Experiment II.

\subsection{Experiment II: Gene expression}\label{experim2}

With this second data example, we want to show that our results make sense at the level of genetically regulated gene expression (i.e., mRNA). For this task, as explained in the Supporting Information, only some variables in our schizophrenia database can be used, comprising some $L=6456$ SNPs that regulate gene expression in the brain, but not in any other tissue of the human body, according to data from the \citet{GTEx}.

We now apply our procedure as in Experiment I, seeing that there are significantly more pairs in putative interaction than in spurious one: 1272 versus 1137 (with \textit{p} of 0.032), after applying the physical distance threshold of 1 Mb. These pairs represent 1539 and 1439 unique SNPs respectively, again a significant difference ($p$-value $\approx0.019$, which drops to 0.0024 by removing SNPs in both sets).

We finally order the $p$-values we obtained for each of the $\comb{L}{2}$ tests we performed on cases, and do the same for controls. We then take the absolute difference of both ranks for each SNP pair. We hypothesize that these absolute rank differences will tend to be greater on the true positive list than on the false positives. We perform a Wilcoxon--Mann--Whitney $U$ test and we find that we can confirm that it is the case, with a $p$-value of less than $2.2\times10^{-16}$.

All in all, the results of Experiment II indicate that we are detecting some genuine signal at the level of genetically regulated gene expression. This would be very unlikely if our method did not function correctly.

\section{Discussion and Conclusion} % Literal name required by the BIMJ.

Distance correlation has been shown to characterize independence in certain metric spaces, establishing how this theory is valid for the case of 3-point marginal spaces. With this approach, a hypothesis test based on the general characterization of independence that distance correlation offers has been designed, extending the idea of LCTs \citep{Cai:Liu} to ternary data.

We derive the explicit asymptotic null distribution of the distance-covariance statistics that arise. To our knowledge, the usage of distance correlation in discrete spaces (in genomics or elsewhere) ---and, in particular, its application to the search for SNP-SNP interactions--- has no precedents in the literature. Moreover, no previously published research has attempted to perform large-scale multiple testing with any of the techniques derived from energy statistics \citep{TEOD}. However, what does exist in the literature is the usage of distance correlation for finding the association between genetic data (as observations of continuous random variables in Euclidean spaces) and a phenotype \citep{Hua}, which is another interesting problem, but completely different both regarding biological and mathematical factors.

Simulations show that the calibration of significance is adequate and that power is considerably high against various alternatives. We also show that we generally outperform \citet{BOOST} ---one of the most popular epistasis detectors--- in the scenarios we have studied.

The schizophrenia database has been analyzed with our methodology, obtaining biologically sound results at the level of genetically regulated gene expression. Some recent studies show evidence of epistasis between regulatory regions of the human genome \citep{Lin,Patel}, which supports our findings.

In order to frame our results, we would like to emphasize that all popular epistasis detectors find a large number of false positives and do not have really high power \citep{Russ}. Therefore, the main limitation of our method (as it is of any other for this task) is that it is very difficult to make any solid discoveries when working with real data. Epistasis detection is an extremely challenging biostatistical problem, in which there is still much progress to be made, given its key role in human complex genetics \citep{Van Steen:Moore}.

Future work may include investigating further distances, like the pseudometrics at the end of Section~\ref{sec:space3}. Another interesting point would be to adapt our techniques to the search for interactions between the mitochondrial (haploid) and nuclear (diploid) genome of an individual, that is, to the study of independence between binary and ternary variables. More generally, a promising line of future research could be the extension of the philosophy of our current paper to the categorical testing for two variables with an arbitrary number of categories, as in \citet{F:categorical}. Finally, one could even aim at testing for independence between SNP data and a phenotype of interest (regardless of it being a continuous, binary or survival-type outcome), which would amount to testing for marginal effects of individual variants, which is a main goal of GWASs \citep{15y}.

\section*{Acknowledgements}

This work has been supported by project \mbox{PID2020-116587GB-I00}, granted by {MICIU/AEI/10.13039/ 501100011033} (Spanish Ministry of Science). F. Castro-Prado's research was supported by the USC Institute of Mathematics (IMAT) and the German Cancer Research Center (DKFZ). The schizophrenia dataset was generated under the support of the Instituto de Salud Carlos III (grant number ISCIII/PI14/01020) to J. Costas, cofounded by the European Regional Development Fund (ERDF). All authors are grateful to the Galician Supercomputing Centre (CESGA) for access to their facilities. F. Castro-Prado is grateful to all the participants and tireless organizers of the ``Archimedes'' Contest (Spanish Ministry of Science) and the ``Eladio Vi\~nuela'' Molecular Biology Summer School (``Men\'endez Pelayo'' International University), since both events were important sources of inspiration. We also thank Dr Dominic Russ (University of Birmingham) for fruitful exchanges of ideas. All authors would like to show appreciation for the time and effort that the editor, associate editor, referees, reproducible research team and production office have devoted to make our article better. Funding for open access charges has been provided by the University of Santiago de Compostela / CISUG.

\section*{Conflicts of interest}
The authors declare no conflicts of interest.

\section*{Data availability statement}
The data that support the findings of this study are openly available in Zenodo at \url{http://doi.org/10.5281/zenodo.20744952}, as well as in the Supporting Information of the article.

\newpage
%\phantom{aaaa} % Adds extra blank page
\end{document}

%\DOIsuffix{bimj.DOIsuffix}
\DOIsuffix{bimj.70150}
\Volume{68}
\Issue{5}
\Year{2026}
\pagespan{1}{}
%\keywords{Association measures; Distance correlation; Epistasis; Genomics of complex diseases; Schizophrenia.\\[1pc]
%}  %%% semicolon and fullpoint added here for keyword style

\title[Supplement to: Testing for epistasis with distance correlation]{Supplement to:\\Testing for genetic interactions in complex disease\\with distance correlation}

%%%% AUTHORS %%%%
\author[Castro-Prado {\it{et al.}}]{Fernando Castro-Prado\orcidlink{0000-0003-3722-2013}~\footnote{Corresponding author e-mail: {\sf{f.castro.prado@usc.es}}}\inst{,1,2,3}}

\author[]{Javier Costas\orcidlink{0000-0003-0306-3990}\inst{2}}

\author[]{Dominic Edelmann\orcidlink{0000-0001-7467-6343}\inst{3}}

\author[]{\\Wenceslao Gonz\'alez-Manteiga\orcidlink{0000-0002-3555-4623}\inst{1}}

\author[]{David R. Penas\orcidlink{0000-0002-7256-3094}\inst{1,4}}

%%%% POSTAL ADDRESSES %%%%
\address[\inst{1}]{Department of Statistics, Faculty of Mathematics, University of Santiago de Compostela, Galicia, Spain.}

\address[\inst{2}]{Psychiatric Genetics Laboratory, Santiago Health Research Institute (IDIS), University Hospital, Santiago de Compostela, Galicia, Spain.}

\address[\inst{3}]{Biostatistics Department, German Cancer Research Center (DKFZ), Heidelberg, Baden-Württemberg, Germany.}

\address[\inst{4}]{Computational Biology Laboratory, Spanish National Research Council (MBG-CSIC), Pontevedra, Galicia, Spain.}

%\vfill

%% maketitle must follow the abstract.
\maketitle                   % Produces the title.

%\vfill

%% If there is not enough space inside the running head
%% for all authors including the title you may provide
%% the leftmark in one of the following three forms:

%% \renewcommand{\leftmark}
%% {First Author: A Short Title}

%% \renewcommand{\leftmark}
%% {First Author and Second Author: A Short Title}

%% \renewcommand{\leftmark}
%% {First Author et al.: A Short Title}

%% \tableofcontents  % Produces the table of contents.

%
%\doublespacing

\vspace*{.1cm}

This supplement contains additional material to the main body of our article.

\vspace*{1cm}

\renewcommand{\thesection}{\Alph{section}}

\newpage

\section{Statistical approaches to epistasis detection}

The recent development of the ``-omic'' disciplines has been parallel to the creation of bioinformatic tools to process the vast amount of data that these experimental sciences produce. The diversity of the available ``-omic'' software is so large that it has even been necessary to develop meta-tools to index the existing techniques. For example, the directory of one of them \citep{OMICtools} contains more than $20\,000$, 900 of which are designed for GWAS data analyses, which in turn contain a subset of 100 that are suitable for epistasis detection.

The existence of such a large number of proposed solutions for such a specific task owes to the surprisingly wide variety of statistical methods that are valid for it---e.g. linear models (standard and generalized), logistic regression, tests on Pearson's correlations, permutation tests, Bayesian inference, random forests, Markov chains, co-information indices, graph theory, and maximal entropy probability models.

Another cause of that diversity of alternatives is the fact that some of the available techniques only focus on a specific subproblem (pairwise gene-gene interactions versus higher orders, binary versus continuous response variable, pedigrees, stratified populations and so forth) and that many different computing strategies can be used to yield results in reasonable amounts of time (for instance, biology-based initial filters, code parallelization, graphical processing units, Boolean operations, machine learning approaches, and ant colony optimization algorithms).

\begin{table}[!htbp]
	%	\vspace{.7cm}
	\centering
	\caption{Some remarkable epistasis detection tools for GWAS data analysis}
	\begin{adjustbox}{width=.95\textwidth}
		%		\fbox{
			\begin{tabular}{@{}llll@{}}
				\toprule
				Tool & Statistical techniques & Computational tricks & Reference \\
				\cmidrule(l{-1pt}r{-1pt}){1-4}
				2S-LRM  & Logistic regression & Pre-filtering & \citet{Vaart} \\
				AntEpiSeeker & $\chi^2$ tests & Ant colony optimization & \citet{AntEpiSeeker} \\
				BEAM  & Bayesian MCMC & None & \citet{BEAM} \\
				BOOST & Logistic regression & Boolean operations, parallelization & \citet{BOOST} \\
				BiForce & Linear regression & Boolean operations, parallelization & \citet{BiForce} \\
				CES   & Evolutionary algorithms & Artificial intelligence & \citet{CES} \\
				CINOEDV & Information theory & Swarm intelligence on hypergraphs & \citet{CINOEDV} \\
				EpiGPU & Linear regression & GPU architectures & \citet{EpiGPU} \\
				EpiACO & Information theory & Ant colony optimization & \citet{EpiACO} \\
				EpiBlaster & Pearson's correlations & GPU architectures & \citet{EpiBlaster} \\
				Fiúncho & Information theory & Parallelization & \citet{Fiuncho} \\
				GLIDE & Linear regression & GPU architectures & \citet{GLIDE} \\
				GWIS  & ROC curve analysis & GPU architectures & \citet{GWIS} \\
				IndOR & Logistic regression & Pre-filtering & \citet{IndOR} \\
				MDR   & Combinatorics, resampling & Pre-filtering & \citet{MDR} \\
				Random Jungle & Random forests & Parallelization & \citet{Random Jungle} \\
				SNPruler & Information theory & Branch and bound algorithms & \citet{SNPruler} \\
				Stage-wise LRT & GLMs, closed testing & Hierarchical testing & \citet{Franberg} \\
				Wtest & Logistic regression & None & \citet{wtest} \\
				\bottomrule
			\end{tabular}
			%		}
	\end{adjustbox}
	\label{HumGenet}
\end{table}

Table~\ref{HumGenet} summarizes some of the existing methods, including the ones reviewed by \citet{Gusareva}, \citet{Niel} and \citet{Russ} and some others that we consider representative. Some are very widely used, like BOOST, due to its being implemented in the popular genetics toolset PLINK \citep{PLINK}; whereas other methods in the table have not been used much in practice.

\section{Derivation of the asymptotic null distribution of the test statistic}\label{sec:proofs}

We first state, as Lemma~\ref{lem:discrete:dist}, a result by \citet{DJ} that is instrumental in the proof of Theorem 3.1, which is the central theoretical result of our article. When it comes to proving Theorem 3.2 (which we do not use in simulations nor in real data analyses), we make use of work by \citet{HH}. We refer the reader to that article ---instead of reproducing its content here---, in order to stay on topic. Finally, we provide the proof of Corollary 3.3., which is the key result for the application of our methodology.

\subsection{A lemma for the discrete metric}\label{subsec:lemma:discrete}

\setcitestyle{square}
\begin{lemma}[\cite{DJ}, Theorem 7] \label{lem:discrete:dist}
Let $X$ and $Y$ be random variables with supports $\{1,2,\ldots,I\}$ and $\{1,2,\ldots,J\}$, respectively; where $I,J\in\Zplus$.

We now construct the entries of matrix $\bL^X=(L^X_{rs})_{r,s}\in\R^{I\times I}$ as follows:
$$
L^X_{rs}=p_s\left( \delta_{rs}-p_r-p_s+\sum_{i=1}^I p_i^2\right) ,
$$
where $\delta_{\cdot\cdot}$ is the Kronecker delta and $p_i:=\Prob(X=i)$ is the probability mass of $X$ in $i\in\{1,\ldots,I\}$. Furthermore, let us denote by
$$\lambda^X_1,\ldots,\lambda^X_{I-1}$$
the nonzero eigenvalues of $\,\bL^X$. If $\{\lambda^Y_j\}_{j=1}^{J-1}$ are defined analogously, then the following limit in distribution holds under independence of $X$ and $Y$, as $n\to\infty$:
$$
n \, \widehat{\dCov}_{\text{discrete}}^2(X,Y) \stackrel{\mathcal{D}}{\longrightarrow}
\sum_{i=1}^{I-1}\sum_{j=1}^{J-1}\lambda_i^X\lambda_j^Y Z_{ij}^2 ;
$$
with $\{Z_{ij}\}_{i,j}$ being IID standard Gaussian.

\end{lemma}
\setcitestyle{round}%

\subsection{Proof of Theorem 3.1 (discrete metric)}
Applying Lemma~\ref{lem:discrete:dist}, one gets that, as $n\to\infty$, 	
$$
n \, \widehat{\dCov}_{\text{discrete}}^2(X,Y) \stackrel{\mathcal{D}}{\longrightarrow}  
\lambda_1 \, \mu_1 Z_{11}^2 + \lambda_1 \mu_2 Z_{12}^2 + \lambda_2 \mu_1 Z_{21}^2 + \lambda_2 \mu_2 Z_{22}^2,
$$
where $\{Z_{ij}^2\}_{i,j}$ are IID $\chi^2_1$; whereas $\{\lambda_j\}_{j=1}^2$ and  $\{\mu_j\}_{j=1}^2$ are the nonzero eigenvalues of certain matrices. Namely, $\lambda_1$ and $\lambda_2$ are the non-null eigenvalues of the following matrix:
$$
\bA=\begin{pmatrix} (1-2 p_0 + \sum p_j^2) \, p_0 & ( -p_0 - p_1 + \sum p_j^2) \, p_1 & ( -p_0 - p_2 + \sum p_j^2) \, p_2 \\ (-p_0 - p_1 + \sum p_j^2) p_0 & (1- 2 p_1 + \sum p_j^2) \, p_1 & ( -p_1 - p_2 + \sum p_j^2) \, p_2 \\ (-p_0 - p_2 + \sum p_j^2) p_0 & (-p_1 - p_2+ \sum p_j^2) \, p_1 & ( 1 - 2 \, p_2 + \sum p_j^2) \, p_2\end{pmatrix}.
$$
By multiplying each of the rows of $\bA$ by $p_0$, $p_1$ and $p_2$ respectively, and adding up the rows, one easily sees that this matrix is singular. Using the relation $p_0+p_1+p_2=1$, we can write the characteristic polynomial of $\bA$ as:
$$
P(\lambda) = \lambda \, \left( \lambda^2 - \left( 1-\sum_{j=1}^3 p_j^2\right)  \lambda + 3 \prod_{j=1}^3 p_j\right), 
$$
where $P(\lambda)$ is defined using the sign convention that makes it monic: $\det(\lambda\,\bI_3-\bA)$.

Calculating the roots of $P(\lambda)$ yields $\lambda_1$ and $\lambda_2$. The derivation of $\mu_1$ and $\mu_2$ follows the same strategy.

\subsection{Proof of Theorem 3.2 (Euclidean distance)}
Throughout the proof, we will use the notation:
$$
M := 2 p_0 p_1 + 2 p_1 p_2 + 4 p_0 p_2 = 2 p_0 (1-p_0) + 2 p_2 (1-p_2).
$$
Applying Theorem 4.12 and Lemma 4.14 in \citet{HH}, we obtain that, when $n\to\infty$,
$$
n \, \dCovh_{\text{Euclidean}}^2(X,Y) \stackrel{\mathcal{D}}{\longrightarrow} \sum_{i,j=1}^3 \lambda_i \mu_j Z_{ij}^2,
$$
as long as $X$ and $Y$ are independent. In the above expression, $\{Z_{ij}^2\}_{i,j}$ are IID $\chi^2_1$ and $\{\lambda_j\}_{j=1}^3$ are the eigenvalues (all of them nonnegative real numbers) of the following matrix:
$$
\begin{pmatrix} (2 p_1 + 4 p_2 - M) \, p_0 & (-1 +  p_1 + 3 p_2 - p_0 +M) \, p_1 & M \, p_2 \\(-1 +  p_1 + 3 p_2 - p_0 +M) p_0 & (2 p_2 + 2 p_0 - M) \, p_1 & ( -1 + 3 p_0 + p_1 + p_2 -M) \, p_2 \\ M p_0 & ( -1 + 3 p_0 + p_1 + p_2 -M) \, p_1 & (2 p_1 + 4 p_0 - M) \, p_2\end{pmatrix}.
$$
On the other hand, $\{\mu_j\}_{j=1}^3$ are defined analogously. Computing the characteristic polynomial and proceeding as we did for Theorem 3.1 completes the current proof.

\subsection{Proof of Corollary 3.3}

As a consequence of (Bernoulli's) strong law of large numbers and the continuous mapping theorem, for each $k \in \{1,2\}$, we have:
\begin{align*}
	\widehat{\lambda_{k}} \stackrel{a. s.}{\longrightarrow} \lambda_{k}, \quad \quad \widehat{\mu_{k}} \stackrel{a. s.}{\longrightarrow} \mu_{k}.
\end{align*}
as $n \to \infty$. Hence, it follows that
$$
\widehat{\lambda_1} \, \widehat{\mu}_1 Z_{11}^2 + \widehat{\lambda_1}  \mu_2 Z_{12}^2 + \widehat{\lambda_2} \widehat{\mu_1} Z_{21}^2 + \widehat{\lambda}_2 \widehat{\mu}_2 Z_{22}^2 \stackrel{a. s.}{\longrightarrow} \lambda_1 \, \mu_1 Z_{11}^2 + \lambda_1 \mu_2 Z_{12}^2 + \lambda_2 \mu_1 Z_{21}^2 + \lambda_2 \mu_2 Z_{22}^2,
$$
as $n \to \infty$. Therefore, for almost all $\mathcal{S}$, when $n \to \infty$,
$$
\widehat{c}_\alpha(\mathcal{S}) \longrightarrow c_\alpha;
$$
with $c_\alpha$ being the $(1-\alpha)-$quantile of the limit distribution. Now, by Slutsky's theorem, under $H_0$,
$$
\frac{n \, \widehat{\dCov}_{\text{discrete}}^2}{\widehat{c}_\alpha} 
\distrilim  \frac{ \lambda_1 \, \mu_1 Z_{11}^2 + \lambda_1 \mu_2 Z_{12}^2 + \lambda_2 \mu_1 Z_{21}^2 + \lambda_2 \mu_2 Z_{22}^2}{c_\alpha}.
$$
This allows us to conclude that, whenever $X$ and $Y$ are independent, for $n \to \infty$,
$$
P\left(\frac{n \, \widehat{\dCov}_{\text{discrete}}^2}{\widehat{c}_\alpha}  \geq  1\right) \longrightarrow  P\left(\frac{ \lambda_1 \, \mu_1 Z_{11}^2 + \lambda_1 \mu_2 Z_{12}^2 + \lambda_2 \mu_1 Z_{21}^2 + \lambda_2 \mu_2 Z_{22}^2}{c_\alpha}  \geq  1\right) = \alpha.
$$\qed

\section{Large-scale correlation tests (LCTs)}

In this section of the supplement, we want to provide a short introduction to the main two philosophies of hypothesis testing for high-dimensional covariance, correlation and precision structures (Section~\ref{2cult}). This serves as an introduction to the LCTs by \citet{Cai:Liu}, a procedure that features prominently in the main article because it is a testing approach similar to ours, but using Pearson's (product-moment) correlation instead of distance correlation. The LCTs are then explained in detail in Sections~\ref{atanh} through~\ref{lctb}.

Let us first introduce some notation. First, given $L\in\mathbb{Z}^{+}$ positions of the genome, let $\mathbf{X}=\left(X_j\right)_{j=1}^L$ and $\mathbf{Y}=\left(Y_j\right)_{j=1}^L$ be the corresponding random vectors containing the genetic information for case and control individuals, respectively. In the setting studied in the main body of our article, each component has support in $\{0,1,2\}$, but for the methods presented in this appendix, based on product-moment covariances, let $\mathbf{X},\mathbf{Y}\in\R^L$.

We denote the correlation matrices for cases and controls as follows.: $$\mathbf{R}_1=(\rho_{ij1})_{i,j}\in\mathbb{R}^{L\times L}\;\;\;\;\text{and}\;\;\;\;
\mathbf{R}_2=(\rho_{ij2})_{i,j}\in\mathbb{R}^{L\times L},$$ 

We observe samples $\{\mathbf{X}_k\}_{k=1}^{n_1}\text{ IID }\mathbf{X}$ and $\{\mathbf{Y}_k\}_{k=1}^{n_2}\text{ IID }\mathbf{Y}$, which are assumed to be independent of each other. As we are in a high-dimensional setting, we assume $L$ to be large, both in absolute terms and when compared with $n_1$ and $n_2$.

\subsection{Testing for high-dimensional covariance, correlation and precision structures}\label{2cult}

As discussed by \citet{AnnuRev:Cai}, there are two main approaches to testing hypotheses on a high-dimensional covariance, correlation or precision matrix. The first one is the \emph{global} testing of a single hypothesis for a very large matrix, originating from random vectors of very high dimension, which introduces evident complications. An example of this would be to test for $$H_0:\mathbf{R}_1-\mathbf{R}_2=\boldsymbol{0},$$
which is of clear interest when one wants to test for differential co-expression. Rejecting such a hypothesis, however, provides no information on which pairs of components (in our setting, SNPs) have differential correlation in the two groups (in our setting, cases and controls).

The other one is the multiple testing of \emph{local} covariance, correlation and network structures. This involves testing hypotheses on pairs of the components of the high-dimensional random vector(s) under consideration, a large number of times (in our setting, $(L^2-L)/2$ times). This comes with statistical and computational challenges that are not present when testing only one of those $(L^2-L)/2$ hypotheses separately. The example that is most interesting to us is the test
\[\begin{cases}{H_{0ij}}:\rho_{ij1}=\rho_{ij2}\\H_{1ij}:\rho_{ij1}\neq\rho_{ij2}\end{cases}\]
for each pair $(i,j)\in([1,L]\cap\mathbb{Z})^2$ so that $i<j$. We will now focus on the LCTs by \citet{Cai:Liu}, which are a robust procedure for testing these hypotheses. The next section motivates the need for transcending conventional correlation testing for this problem.

\subsection{LCT: classical approach}\label{atanh}

A scarcely innovative approach would be to stabilize the variance of the sample correlation coefficients via Fisher's $Z$ transformation (\textit{atanh}). One could think of combining this strategy with a procedure that controls the false discovery rate (FDR), such as the ones by \citet{BH} or \citet{BY}, thus establishing the desired large-scale correlation test (LCT). The main drawback of this idea is that, when normality is not ensured, the behavior of the test statistic differs from the well-known asymptotic distribution of the Gaussian case. Simulation studies \citep{Cai:Liu} show that this method performs very poorly (both with Benjamini--Hochberg and Benjamini--Yekutieli), especially when compared with the LCTs that we introduce next.

\subsection{LCT with normal approximation (LCT-N)}\label{lctn}
The first test that \citet{Cai:Liu} devised, the LCT-N, is based on the test statistic
\[T_{ij}:=\frac{\hat\rho_{ij1}-\hat\rho_{ij2}}
{\sqrt{\frac{\hat\kappa_1}{3n_1}\left(1-\tilde\rho_{ij}^2\right)^2+\frac{\hat\kappa_2}{3n_2}\left(1-\tilde\rho_{ij}^2\right)^2}}\text{,}\]
where $\hat\kappa_1$ and $\hat\kappa_2$ are the respective sample kurtoses of $\mathbf{X}$ and $\mathbf{Y}$, and $\tilde\rho_{ijl}$ is a thresholded version of $\hat\rho_{ijl}$, for $l\in\{1,2\}$; with $\tilde{\rho}_{ij}^2:=\max\{\tilde{\rho}_{ij1}^2,\tilde{\rho}_{ij2}^2\}$.

$H_{0ij}$ will be rejected when $|T_{ij}|$ is greater than a certain threshold $\hat t_{\alpha}\in\Rplus$, which depends on the nominal value $\alpha\in]0,1[$ under which one wants to maintain the FDR. The formula for computing $\hat t_\alpha$ works under the assumption that the initial distributions are Gaussian, or are not very far away from being so (elliptical contours). Consequently, the LCT-N should not be used in other contexts.

\subsection{LCT with bootstrap (LCT-B)}\label{lctb}
If the distributions of $\mathbf{X}$ and $\mathbf{Y}$ are totally unknown, it is reasonable to use resampling techniques in order to approximate the tail of the distribution of $T_{ij}$, which determines $\hat t_{\alpha}$. The bootstrap scheme that \citet{Cai:Liu} built to this purpose is consistent and leads to a threshold $\hat t_{\alpha}^*$, which defines the LCT-B. This test is supported by strong theoretical results, which were proven in the original \citeyear{Cai:Liu} article.

\section{A resampling approach to the testing procedure}\label{resam:appr}

We initially attempted to approximate the $p$-values of our test by using a permutation-based approach, which is the gold standard in the distance covariance literature \citep{TEOD}. We briefly discuss this brute-force approach, as a \emph{negative result} that can be of utility to other scientists. Although theoretically sound, this strategy leads to such high computation times that it is not suitable for most applications. This is particularly true for genomics, which is so data-intensive that not even the high-performance computing (HPC) approach described in the next section makes resampling feasible. Instead, we recommend resorting to asymptotic results, like the ones presented in our article. Nonetheless, in other scenarios (of lower dimensionality, or in which the theoretical derivation of the asymptotic null distribution is not possible), the following approach might be of interest.

In order to approximate the null distribution of the test statistic $\widehat{\dcov}(Z_i,Z_j)$, it is possible to devise a resampling scheme according to the information that is available under the null hypothesis, which in this case is the independence of $Z_i$ and $Z_j$. As a result, the reasonable thing to do is not to resample from $\{(Z_{i,k},Z_{j,k})\}_k$, but to do it separately from $\mathcal Z_i:=\{Z_{i,k}\}_k$ and $\mathcal Z_j:=\{Z_{j,k}\}_k$ (permutation tests). Thus, it suffices to compute $B\in\mathbb{Z}^{+}$ statistics of the form $$\widehat{\dcov}(\mathcal Z_i^{*(b)},\mathcal Z_j^{*(b)})$$ to obtain a Monte--Carlo approximation of the sampling distribution of the empirical distance covariance under $H_{0ij}$. For a faster algorithm, it suffices to permute only one of the marginals: $$\widehat{\dcov}(\mathcal Z_i^{*(b)},\mathcal Z_j).$$

The usage of permutation tests in this context of metric spaces was inspired by the excellent performance of the same scheme in Euclidean spaces \citep{SRB,TEOD,Domi:surv}. It has the drawback that there is not the same kind of fully-fledged formal justification of consistency as the one for the naive bootstrap (see\setcitestyle{square}%
\citet{Arcones:Gine} and page 100 of \citet{Jakobsen}),\setcitestyle{round}%
which should not be a source of concern in practice, like in the Euclidean case. Some preliminary experimental checks (data not shown) point to a better performance of permutation tests versus the naive bootstrap approach in the ternary variable setting, but further studies would be needed to verify it.

It should also be clarified that authors such as \citet{Cai:Liu} and \citet{SRB} argue that the number of resamples $B$ is relatively unimportant for their methods to work, as long as it is not extremely small. With this in mind, and also taking into account that the running time is $O(B)$, we decided to use a moderate value for $B$ in the present article, namely the one proposed by \citet{SRB} as a function of sample size $n$:
\[B(n)=200+\lfloor{5000}/{n}\rfloor\text{,}\]
where $\lfloor\cdot\rfloor:\;\mathbb{R}\to\mathbb{Z}$ is the floor function. Some empirical checks confirm that increasing~$B$ with respect to the value above causes barely noticeable improvements (if any), both in terms of the calibration of significance levels (as long as the nominal value is not extremely small) and of power.

Nevertheless, if the goal was to control the FDR, the situation would be quite different --- both \citet{BH} and \citet{BY} are based on ordering \textit{p}-values and, consequently, $B$ cannot be too small with respect to the number of hypotheses that are tested, in order to avoid an excessive build-up of null \textit{p}-values. Unfortunately, when it comes to analyzing real data, this becomes unfeasible (see the empirical issues in the upcoming Appendix~\ref{comp:chall}). Instead, we opted for controlling the FWER (and not the FDR), which has the slight drawback of neglecting the weak dependencies among the $\binom{L}{2}$ hypothesis tests (both due to each individual SNP being present in multiple pairs and due to linkage disequilibrium).

\section{Computational challenge of the resampling approach}\label{comp:chall}

The implementation of the brute-force testing procedure, as presented in Appendix~\ref{resam:appr}, was an extremely challenging issue from the computational point of view, given the high dimensionality of the data and the large number of hypothesis tests resulting from the combinatorial explosion. Thus, a quite sophisticated set of computer technologies and strategies was required to obtain results in somewhat manageable computational times. 

\begin{table}[!htbp]%[htbp!]%[H]%[!htbp] % `float` package
	\centering
	
	\caption{Comparison of running times for the different versions of the code, all of them referring to the permutation testing approach. It should be noted that the times for runs of the sequential code on the largest GWAS are estimations.}
	
	%	\vspace*{-.2cm}
	%
	%		\fbox{
		\begin{tabular}{@{}lrrr@{}}
			\cmidrule(l{-7pt}r{-7pt}){1-4}
			& Simulation, & GWAS, & GWAS,\\
			Code version & $R=10^3$ & $L=1000$ & $L=4000$\\
			\cmidrule(l{-7pt}r{-7pt}){1-4}
			R sequential & 12 h 10 min & 42 days 1 h & 2 years \\
			R \& C sequential & 3 h 59 min & 2 days 1 h & 30 days  \\
			R \& C parallel & 50 min & 2 h 41 min & 2 days  \\
			\cmidrule(l{-7pt}r{-7pt}){1-4}
		\end{tabular}
		%		}
	
	\vspace*{.2cm}
	\label{Tempos}
\end{table}

As a general rule, any statistical technique based on GWAS data will suffer from the issues that are inherent to such input (high dimension and low sample size). To illustrate this point, Table~\ref{Tempos} compares the run times of the original \emph{R} code with another one, whose core is implemented in the compiled language \emph{C}, this way making the numerical crunching far swifter. This second code ---labeled ``R \& C'' in the table--- also includes some high-performance computing (HPC) improvements and, what is more, it can run in sequential or in parallel mode (i.e., the workload can be distributed among different processors, increasing the computational speed by a factor that is approximately equal to the number of available processors). For a comparison of performance like the one in Table~\ref{Tempos}, it is crucial to carry on the experiments in the same environment, which in our case was supercomputer \textit{Finisterrae II} (Galician Supercomputing Center, CESGA).

Hence, in light of the order of magnitude of these times (the \emph{R} version would need up to two years in large-scale settings, while the \emph{R} \& \emph{C} parallel implementation only requires ten hours), it is fully justified to resort to HPC strategies in a compiled language, especially if one takes into account that a GWAS can easily involve millions of SNPs, with the running time being a linear and monotonically increasing function of $\binom{L}{2}$ and, consequently, $O(L^2)$; as illustrated by the ratio between the GWAS columns of Table~\ref{Tempos}.

For the parallel version of the \emph{R \& C} code, in each case, the lowest amount of hardware that yielded results within a reasonable amount of time was used --- 12 cores for simulations, and 48 processors for real data analyses. To reduce the times by a factor of $f$, it would suffice to increase the number of processors $f$ times, as long as logistic constraints make it possible.

Moreover, the algorithm was parallelized in two alternative ways:
\begin{enumerate}
	\item using a shared memory paradigm via the OpenMP library, distributing the computational effort among the different cores that exist within a processor;
	\item applying a distributed memory strategy, where different computational nodes ---that belong to various machines--- are able to share workload via a message protocol, which in our case is the MPI library.
\end{enumerate}

The first parallelization (that is very easy to implement in the main loop of the algorithm) was useful for applying the test to simulated data, where the dimensionality was not too problematic. However, the number of parallel running threads one can add is limited by the number of available cores in a CPU processor chip, and it is not enough to address real data examples. For this reason, a distributed memory parallelization was developed, with a classical leader/follower paradigm, where hundreds of processors can work together to reduce complexity. It consists in:
\begin{enumerate}
	\item A processor (\emph{leader}) calls \emph{R} routines that load the matrices that contain the input, split it and distribute it among several processors (\emph{followers}).
	\item Each processor works with one fragment of the matrix, running the iterations that have been assigned to it (i.e., performing independence tests for a fraction of the total SNP pairs).
	\item Once each follower finishes its part, it sends the results to the leader.
	\item The leader builds the final \textit{p}-value matrix, which is later used to wrap up the results in \emph{R}.
\end{enumerate}

Please note that ``{leader/follower}'' is just one of several proposed replacements for the historical name of this technology, which is linked to racial discrimination and we therefore avoid reproducing such nomenclature here.

The \emph{R \& C} version of our code combines an interface in the programming language \emph{R} with a core in \emph{C}, with the latter being devoted to performing low-level computations in a time-efficient manner. Another important factor that decreases running times in our implementation is the use of specific libraries to codify low-level operations that involve large vectors and matrices. Namely, the well-known Intel MKL libraries and SIMD (Single Instruction; Multiple Data) techniques have been applied to exploit data-level parallelism --- using an extension in the registers and the arithmetic and logic instructions present in modern microprocessors, they can process the same operation simultaneously on the elements of an array through a single instruction. We found it particularly useful when it came to implementing matrix operations.

It was not possible to resort to preexisting software because the most efficient distance-correlation-related algorithms \citep[like the one by][]{Chaudhuri} are only designed for the Euclidean case and, therefore, not adaptable to the structure of the 3-point spaces that are the scope of our article.

\section{Genomic database}
The SNP data used in our article come from a case-control study of schizophrenia, which was performed on $n_1=585$ patients and $n_2=573$ control blood donors, all of them of Galician origin, as described by \citet{Galicia}.

\begin{sloppypar}
Each individual's genome was sequenced using microarray \mbox{\textit{PsychArray-24 BeadChip}} (Illumina, San Diego, California). After genotyping, several conventional quality controls were performed. Namely, to avoid experiment-derived problems, we excluded from the database every SNP that verified at least one of the following conditions:
\end{sloppypar}
\begin{enumerate}
	\item The minor allele frequency (MAF) is less than 1 \% in our samples.
	\item The genotype proportions differ significantly from Hardy--Weinberg equilibrium in the control sample, for nominal $\alpha=0.001$.
	\item The \emph{call rate} (proportion of non-missing data after genotyping) is under 95 \%, or it is significantly different between cases and controls ($p-$value of less than 0.001).
\end{enumerate}

The previous filters aim at only retaining those SNPs for which one can clearly discriminate between the three possible genotypes.

Individuals were removed when, after the SNP quality control, the genotyping for more than 5 \% of their SNPs was missing. For assessing cryptic relatedness, we computed the identity by descent proportion (pi-hat statistic) for each pair of individuals and, for every pair in which $\hat\pi>0.15$, one of its members was removed.

All the aforementioned restrictions were applied with the default algorithms and implementations for GWAS quality controls of PLINK \citep{PLINK}.

The final SNP count for the schizophrenia database, which is used in Section 5 of the main article for illustrating our methodology, is $L=6,371,078$. It contains observations of $n_1=585$ cases and $n_2=573$ controls.

Data collection followed the guidelines of the Declaration of Helsinki, was approved by the Galician Ethical Committee for Clinical Research, and
participants signed an informed consent; as stated in \citet{Galicia}.

\section{Technical details of the real data application}

We now explain some non-essential technicalities that were left out of the explanation of the experimental procedure in Section 5 of our article.

\subsection{Experiment I: Functional enrichment}

For this experiment, among all the autosomal SNPs, only missense SNPs (i.e., those that induce a change in the amino acid sequence of a protein) were considered, as a way to detect interaction at the protein level.

To determine which SNPs are missense and which not, our reference was the  ENSEMBL Biomart database \citep{Biomart}. Since our schizophrenia data refer to the GRCh37.p13 (hg19) assembly of the human genome, we chose the Biomart version accordingly. For a comprehensive review of ENSEMBL and Biomart, we refer the reader to \citet{ENSEMBL:rev}.

As we wanted to have the ability to detect some signal, and it is remarkably difficult to detect any instance of epistasis in real data \citep{Russ}, we restricted ourselves to those SNPs for which the least frequent of the two alleles is observed with a frequency of at least $0.05$ in our samples and $0.01$ at ENSEMBL Biomart. The latter filter also ensured that the missense SNPs under study are well known and annotated. This yielded a SNP count of $19,356$.

We then removed all the SNPs that lie in any of the 25 small regions of the human genome which are known to suffer from \emph{long-range linkage disequilibrium} (LD), a phenomenon that can cause confusion between true SNP-SNP interactions related to schizophrenia and associations due to the architecture of chromosomes \citep{HighLD:ajhg}. The list of those 25 regions can be found, for example, in \citet{Facal}.

Furthermore, once the high LD regions were removed, we followed standard practice among geneticists \citep{Abdellaoui:2013} and, from the remaining $18,268$ SNPs, we \emph{pruned} those showing evidence of short-range LD, setting $r^2<0.1$ in PLINK 1.9 \citep{PLINK}. We used windows of width 500 SNPs, shifting them $+1$ SNP at each step. Thus, we obtained the $L=8030$ SNPs that were analyzed with our method in Experiment I of the main body of the article.

Whenever we had to assign SNPs to genes, we once again followed the annotations at ENSEMBL Biomart.

\subsection{Experiment II: Gene expression}

In Experiment II, we only used SNPs that regulate gene expression (i.e., eQTLs) in any of the brain tissues, and nowhere else in the human body, according to data from the \citet{GTEx}. With this aim, we downloaded \texttt{GTEx\textunderscore Analysis\textunderscore v7\textunderscore}\texttt{eQTL.tar.gz} (single tissue cis-eQTL data for GTEx Analysis V7, dbGAP accession phs000424.v7.p2) from

\begin{center}
\url{https://www.gtexportal.org/home/datasets},
\end{center}
which again refers to the GRCh37.p13 (hg19) assembly of the human genome. There we simply removed anything that is not a SNP (e.g., insertions), created lists of SNPs that regulate gene expression in brain and non-brain tissues, and performed a set difference.

We therefore chose those SNPs present in our study which are also among the $97,913$ SNPs that act as eQTLs only in brain (and not in other tissues) according to GTEx. Removing data from sexual chromosomes, as well as the high LD regions (as in Experiment I), the SNP count dropped to $56,395$. After once again pruning the short-range LD, we got the final SNP list for this experiment, of length $L=6456$.